\documentclass{birkjour}

 \theoremstyle{definition}
 
 \theoremstyle{remark}

 \numberwithin{equation}{section}

\begin{document}

\title[]
 {Klein, Lie, and their early \\Work on  Quartic Surfaces}

\author[]{David E. Rowe\\Mainz University}




\maketitle




\section*{Introduction}

Special types of quartic surfaces were much studied objects starting in the 1860s, so it should come as no suprise that these played a major role in the early work of Sophus Lie and Felix Klein. This topic was very much in the air when they first met each other in Berlin and attended the famous seminar led by Ernst Eduard Kummer and Karl Weierstrass. Only a few years before, Kummer had published a series of important results on quartic surfaces, including the special ones that today bear his name. Somewhat ironically, it was he and not his former colleague, Jakob Steiner, who first brought Steiner surfaces to the attention of the mathematical world. In what follows, I will begin by briefly recounting this story as well as Lie's early interest in Steiner surfaces, as revealed in his correspondence with Alfred Clebsch. 

Klein's background as an expert in  line geometry led him to investigate the role of Kummer surfaces in the theory of quadratic complexes, but also the closely related complex surfaces first studied by his former teacher, Julius Pl\"ucker \cite{Clebsch 1872}. In G\"ottingen, Klein befriended Max Noether and learned from him about a technique for mapping the lines of a linear complex to points in 3-space, a mapping closely related to another discovered independently by Lie. Thus, once Klein and Lie met in the fall of 1869, their common interests in these geometrical ideas gave them an enormous impulse to pursue them further. These interests eventually led to Lie's theory of transformation groups, as well as to the  new vision for geometrical research that Klein outlined in his ``Erlanger Programm'' (\cite{Hawkins 1989}, \cite{Rowe 1989}).   These famous developments have received considerable attention in the secondary literature, whereas the background story, which concerns  concrete problems in line geometry, has been largely overlooked. This essay aims to bring some of these key background  events to light by focusing
 on the role of  special quartic surfaces in the context of the brief period during which Klein and Lie began their fruitful collaboration.

\section*{Kummer's Rediscovery of Steiner Surfaces}

Jakob Steiner approached the study of algebraic surfaces using purely synthetic methods, in sharp contrast to other leading contemporary geometers, in particular Pl\"ucker and Otto Hesse.
 In the case of cubic surfaces, Steiner devised a number of ways to construct these by starting with the 27 lines that lie on them. One such construction, presented in \cite{Steiner 1857}, begins  with two pairs of trihedral planes that meet in 9 lines. Taking planes through a fixed point $P$, Steiner notes that each contains a cubic curve that passes through $P$ and the 9 points of intersection with the 9 fixed lines. These planes then sweep out a cubic surface  containing these lines. Steiner then obtains the other 18 lines by noting that among the first 9 there will be 6 triples of mutually skew lines. Each such triple then determines a system of generators for a hyperboloid, which will have three additional lines in common with the cubic surface. This, then, accounts for all 27 lines, a construction he presented  without any hint of a proof.

When it came to the famous quartics that bear his name, Steiner was even less forthcoming. His inspiration was again a purely synthetic construction, an approach closely related to the one later taken by Lie. Since plane sections of quartic surfaces yield quartic curves, the question he asked was whether such a surface could contain infinitely many conics. Clearly, these would then appear in pairs representing decomposable quartic curves. Evidently Steiner first imagined such an object in 1844, when he was vacationing  in Rome with Jacobi,   Dirichlet, and the latter's wife Rebecka, sister of the composer Felix Mendelssohn-Bartholdy (see \cite{Rowe 2017b}, introduction to Part I). Jacobi had only recently left K\"onigsberg to accept an appointment as a salaried member of the Prussian Academy, but he was first allowed to spend some weeks in Italy to restore his health.
It seems Jacobi and Steiner got along very well, as they often spent  evenings together talking about 
 various geometrical problems. Probably they had  difficulty 
 understanding each other, given their sharply contrasting mental abilities. Steiner had 
 scarcely any algebraic or analytical skills, whereas Jacobi was, of course, a virtuoso when it came to wielding such instruments. It seems likely that Steiner would have told Jacobi about the exotic quartic, today  often called the Roman surface, that he discovered at this time, but we shall probably never know. In fact, Steiner never published anything at all about this discovery, which only came to light two decades later.

Just a few months after Steiner's death on April 1, 1863, Kummer presented a paper to the Berlin Academy \cite{Kummer 1863},  containing several new results concerning quartic surfaces with families of conics. Among these surfaces was a quartic whose tangent planes cut the surface in pairs of conics. This  turned out to be a rational surface containing three double lines that meet in a triple point. 
Kummer apparently had no idea that he had rediscovered the Steiner surface, but while  writing up his results his colleague, Karl Weierstrass, informed him that Steiner had stumbled upon this special quartic surface many years earlier. Already as a Gymnasium student in Paderborn, 
 Weierstrass had studied Steiner's  early works  in Crelle's Journal. So even though, like Jacobi, his own strengths lay in analysis, Weierstrass had a deep appreciation for Steiner's achievements as a geometer. 
Both he and Kummer also  knew how eccentric this self-taught Swiss genius had been, so perhaps it came as no surprise to them that he never published his discovery;  in fact, he nearly took this secret with him when he went to his grave.

Kummer's paper was submitted to the Academy on July 16, but Weierstrass added a short  note  explaining the  idea behind Steiner's purely synthetic construction \cite{Weierstrass 1863}. He  
wrote this up based on his recollections of a conversation with Steiner that took place perhaps a year or two before the latter's death, noting that  he thought it  unlikely  Steiner would have left behind any written account of his discovery.
At the time of Steiner's death, Weierstrass was nearly the only one in Berlin who enjoyed friendly relations with him, as he had become an increasingly embittered old man who likened himself with a ``burnt-out volcano'' (\cite{Biermann 1988}, 84).
Weierstrass even  promised the feisty Steiner that he would continue teaching  courses in synthetic geometry once his elderly colleague could no longer do so. He kept that promise; during the decade that followed Steiner's death Weierstrass taught courses in geometry no fewer than seven times.

\begin{figure}[h]
        \centering 
        \includegraphics[width=7.5cm,height=6.5cm]{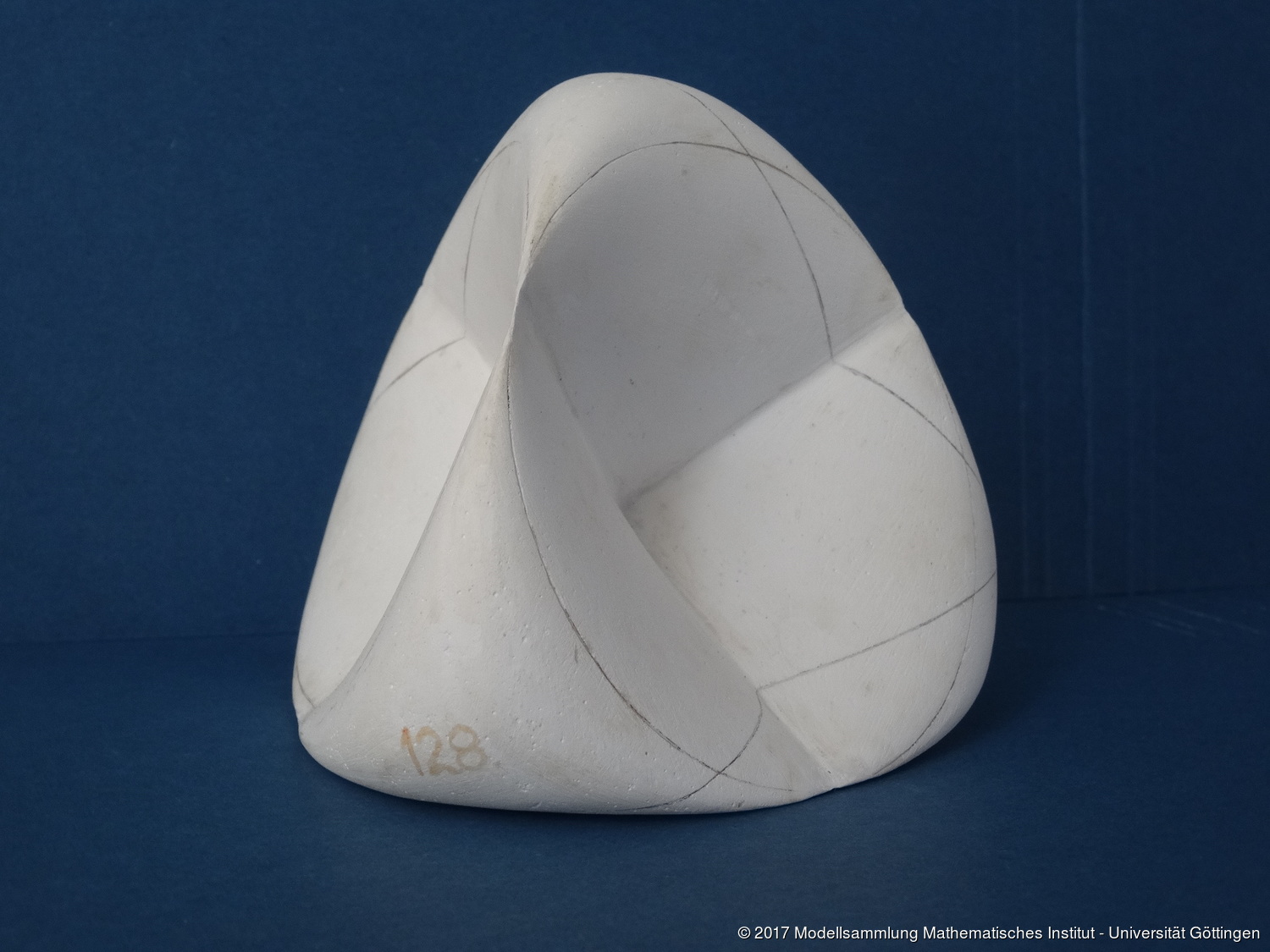}
        \caption{Kummer's model of a Steiner surface (Courtesy of the Collection of Mathematical Models,  G\"ottingen University)}
\label{fig:maya}
\end{figure}

Later that same year, Kummer submitted a short note to the Berlin Academy on the Steiner surface, accompanied by a plaster model illustrating the case:  $$y^2z^2+z^2x^2+x^2y^2-2cxyz=0.$$
 Somewhat ironically, his discovery of the Steiner surface actually preceded his more celebrated work on so-called Kummer surfaces -- quartics with 16 double points and 16 double planes -- which began shortly after this (\cite{Kummer 1975}, 418-432). In fact, all of  his work on quartic surfaces during the 1860s was closely tied to his studies of ray systems in geometrical optics, work inspired by the earlier researches of William Rowan Hamilton. This included the special quartic known as 
 the Fresnel wave surface, which  Kummer himself identified as a ``Kummer surface,'' but with only four real double points (the remaining 12 being imaginary). 

As for the ``Steiner surface,'' Kummer's simple analytical presentation of it led to a flurry of subsequent investigations  by Schroeter, Clebsch, Cremona, Laguerre, and in 1878, Sophus Lie.
After two decades, Steiner's name was securely attached to this object, even though Weierstrass's memory was the only known link connecting Kummer's work with what Steiner had earlier conceived during his stay in Rome.
Then, in 
 the early 1880s, the  inspiration behind Steiner's discovery finally came to light. Weierstrass, who was the editor of Steiner's Collected Works, at some point  found  a brief note in the geometer's papers that described the original construction. Except for a few details, Weierstrass now realized that his recollection from 1863 had been quite accurate. So in 1882 he decided to include this sketchy note as an appendix to the second volume of Steiner's {\it Werke} (\cite{Steiner 1882}, 723-724), this time adding  a lengthier recollection of his discussions with Steiner about this surface (\cite{Steiner 1882}, 741-742). 
Here, for the first time, we read when and where Steiner  originally found this surface, namely in 1844 while in Rome; indeed, we learn that it was for this reason that he liked to call this quartic his ``Roman surface,'' the other name by which it is known today.
Weierstrass further relates that Steiner was truly tortured by the question of whether his construction actually led to a quartic surface. He suspected, in fact, that his  Roman surface might  contain an additional  imaginary component as a kind of  ``ghost image.''  So  he asked  Weierstrass  to clarify this matter analytically, and the latter was able to do so easily by giving  a rational parameterization of the surface using  homogeneous quadratic polynomials.

\section*{Lie's Theorem on Steiner Surfaces}

In the year of Steiner's death, Otto Hesse  published a short obituary, honoring  Steiner as the ``leading geometer of his day''  (\cite{Hesse 1863}, 199).  Referring to Steiner's later period, Hesse  saw this as marked by his struggle with the imaginary, or as Steiner liked to say, those ``ghosts'' that hide their truths in a strange geometrical netherworld. As Hesse clearly saw, this was a fight he could not win, especially given the constraints of his purely synthetic approach to geometry. Still, he was clearly awed by Steiner's daring vision, likening his numerous unproved results with the many puzzles Fermat bequeathed to the mathematical world. 

 Sophus Lie was another geometer who struggled to visualize objects defined by equations in two or more complex variables (\cite{Rowe 1989}, 224-226). He also took an early interest in Steiner surfaces and found a remarkable theorem about them that  drew on Steiner's construction as given in \cite{Weierstrass 1863} as well as analytical arguments presented in \cite{Clebsch 1867}. Although he had already found this result  in 1869, Lie waited almost a decade before publishing it in \cite{Lie 1878}, after which his discovery 
 led to several subsequent investigations by Italian geometers (see \cite{Lie 1934}, 792-793).
 Lie's theorem states that any given  Steiner surface will determine a 3-parameter family of such surfaces, and that one can obtain these  merely by taking arbitrary fixed planes $\Pi$  in space. Lie observed that each tangent plane $T_P$ touching a Steiner surface $S$ cuts out a pair of conics, and these give rise to two pole points in $T_P$: these being the poles corresponding to the polar line $\ell_P = \Pi \cap T_P$  with respect to the two conics. Lie then showed that by varying the planes $T_P$ over $S$, these pole points sweep out a new Steiner surface. The only exceptional cases are those where the plane $\Pi$ is itself tangent to $S$, in which case the construction leads to a quadric surface. Combining these two cases, then, Lie's theorem says that every Steiner surface $S$ determines $\infty ^3$ Steiner surfaces and $\infty ^2$ quadric surfaces.

Lie added a footnote to his paper, in which he pointed out that he had already communicated these results to the University of Christiania in 1869. However, Friedrich Engel and Poul Heegaard, the editors of the first volume of his {\it Gesammelte Abhandlungen},  were unable to find any trace of such a communication (\cite{Lie 1934}, 790). Still,  indirect confirmation that Lie was in possession of this result  back in 1869 can be found in a letter he received in that year from Alfred Clebsch, a document that can be found in Lie's Nachlass in Oslo (see the transcription in the appendix). Lie had  written to Clebsch, asking him if he was aware of this theorem; he also asked Clebsch  about a second result concerning Pl\"ucker's special quartic surfaces that had left him puzzled. Clebsch wrote back to inform him that he had not known Lie's theorem; however, he was able to confirm its correctness  by utilizing a parametrization of the given Steiner surface. This reply to Lie's queries was written, in fact,  just as Lie was about to meet Felix Klein for the first time in Berlin. Knowing this, Clebsch suggested that Lie  speak with  Klein about  the puzzling remarks in  Pl\"ucker's book \cite {Pluecker 1868}, about which Lie had written. 

\section*{Kummer Surfaces}

Felix  Klein learned about various types of special quartic surfaces as a student of Julius Pl\"ucker in the mid-1860s. This was just around the time that E. E. Kummer began publishing on the  quartic surfaces that now bear his name. These Kummer quartics have 16 singular points that lie in groups of six in 16 singular planes, or tropes. Each of these 16 tropes is a tangent plane to the surface that touches it along a conic section containing six of the 16 singularities. These sets of six points can thus be seen to lie in a special position, since only five coplanar points in general position determine a conic. In fact, the  singular points and planes of a Kummer surface   form a symmetric (16,6) configuration, which means that six singular planes also pass through each of the singular points of the surface.

\begin{figure}[h]
        \centering 
        \includegraphics[width=8.5cm,height=7.0cm]{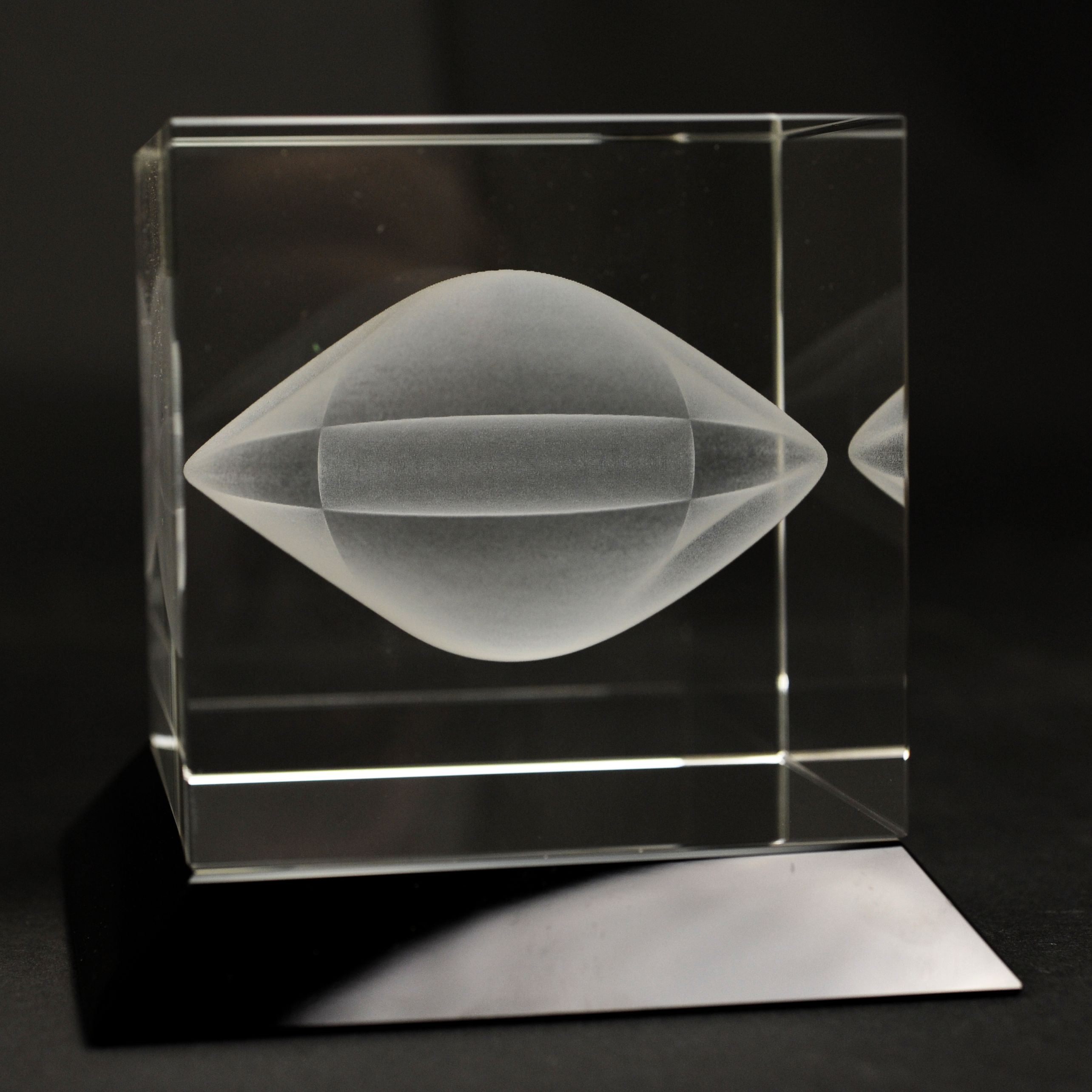}
        \caption{Model of a Fresnel wave surface using laser-in-glass. Courtesy of Oliver Labs.}
\label{fig:maya}
\end{figure}

As noted above, the first such quartic to be  studied in considerable detail was the Fresnel wave surface, which 
 happens to be connected with a famous problem in optics: the phenomenon of double refraction. This arises when a light beam entering a crystal  is refracted in two different directions, so that it splits into two beams in the course of passing through the glass.  In the 1820s Fresnel derived an equation for the wave front of light passing through biaxial crystals; this turned out to be a two-leaved surface described by 
 a quartic equation, though many of  its  deeper properties alluded him. James MacCullagh and 
W. R. Hamilton were  the first to recognize the significance of this structure for optics \cite{MacCullagh 1830}. Hamilton saw that the Fresnel surface has four singular tangential planes that correspond  to its four singular points. A singular plane touches the surface all along a curve, which for quartics will typically be a conic section.
For the Fresnel surface, the singular planes touch it along four circles that lie in two pairs of planes parallel.
 Soon afterward, the 
 Dublin geometer 
George Salmon  
 showed that the Fresnel surface was of the fourth order and class, owing to the fact that it has 16 double points and 16 double planes. The order of a surface is given by its degree in point coordinates, so a quartic has order $n=4$. Geometrically this means that an arbitrary line  meets the surface in four points. Similarly, the class $k$ of a surface is the number of tangent planes through a generic line. In general, the class $k$ of a surface of order $n$ is given by $k=n(n-1)^2-2d$, where $d$ is the number of double points \cite{Salmon 1847}. Hence, for $n=4$ and $d=16$, we get $k=4(3)^2-2(16)=4$, which shows that $d=16$ is the maximum number of double points possible for a quartic surface.

\begin{figure}[h]
        \centering 
        \includegraphics[width=10cm,height=7.5cm]{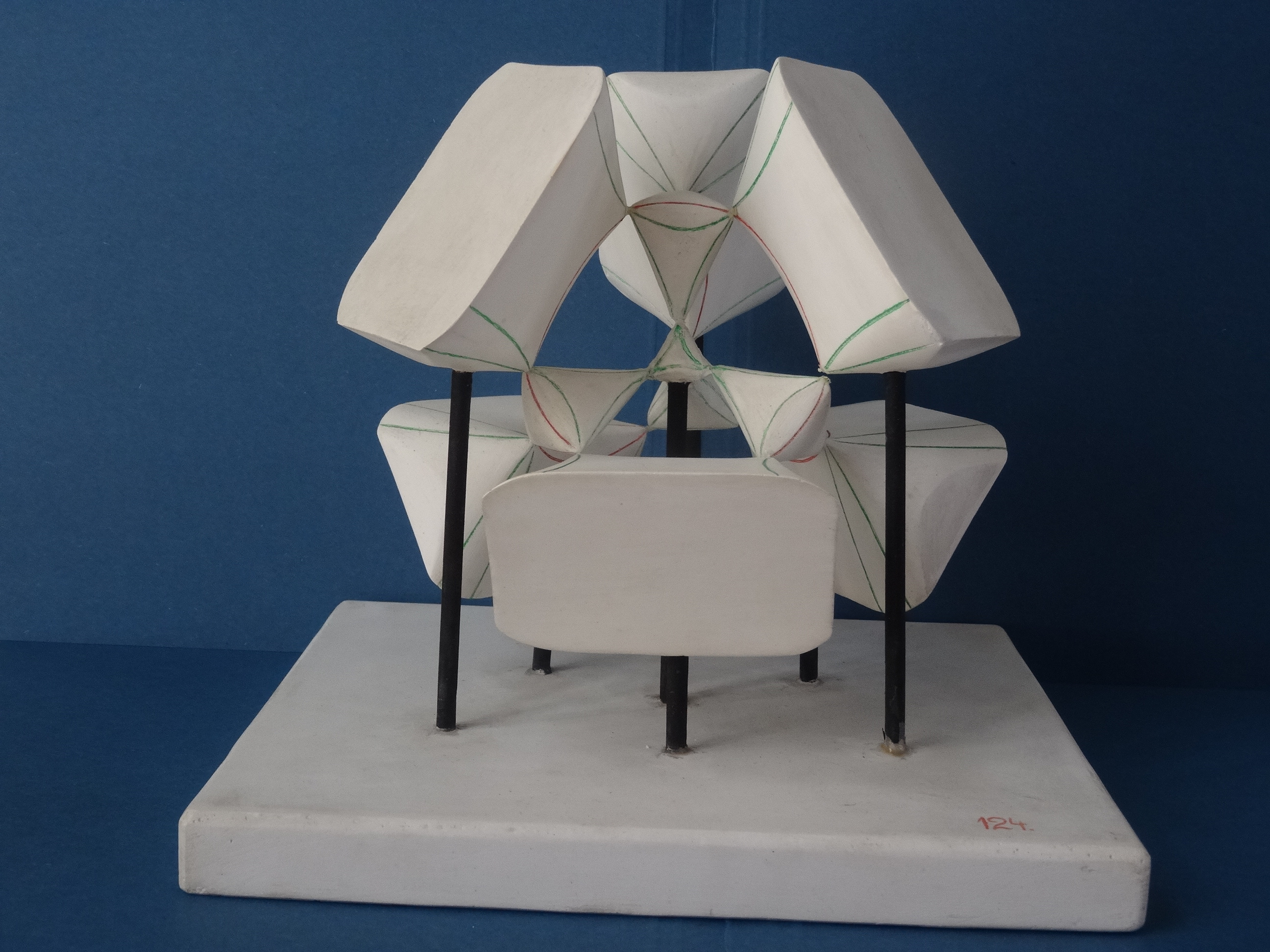}
        \caption{Karl Rohn's model of a Kummer surface with 16 real nodes. Courtesy of the Collection of Mathematical Models,  G\"ottingen University.}
\label{fig:maya}
\end{figure}

Klein was the first to recognize the close connection between Kummer surfaces and another special type of quartic, the complex surfaces  of Pl\"ucker (see below and \cite{Klein1870},  \cite{Klein1874}). In 1871, he designed four zinc models that illustrate the relations between these two types (\cite{Klein1921-23}, II: 7-10). His model of a Kummer surface was similar to the plaster model shown here, which was  built by his student Karl Rohn \cite{Rohn 1877}. 
This highly symmetric quartic can be visualized as consisting of eight tetrahedral-like pieces: an inner tetrahedron with four others attached to each vertex, plus three outer tetrahedra that join at infinity. For the model, these last three tetrahedra have been split and thus truncated  into two pieces by planar sections, which is why the model has six outer pieces. By imagining the opposite pieces to extend through infinity, they would then join to form the three outer tetrahedra (see \cite{Rowe 2018}).

\section*{Quartic Surfaces in Line Geometry}

During the last years of his life, 
Pl\"ucker  devoted a great deal of attention to studying and building special quartics that he called complex surfaces (see \cite {Pluecker 1868}, \cite {Pluecker 1869}). This name comes from line geometry, in particular the theory of quadratic line complexes $K_2$ .  Because these 3-parameter families of lines are very difficult to visualize in toto, 
Pl\"ucker 
studied their local structure instead. His idea was to  fix a  line $g\notin K_2$ and consider all the lines in $K_2$ that meet $g$. This amounts to forming 
$K_2  \cap K_1(g)$, where $ K_1(g)$ is the first-degree complex consisting of the lines in space that intersect $g$. The intersection of two algebraic  line complexes yields a 2-parameter family of lines, called a line congruence. In this case the congruence is of the second order and class since two lines lie in an arbitrary plane and two pass through a generic point. 

Such systems were familiar from earlier work in ray  optics, the background for Kummer's work in the 1860s. Thus it was known that a congruence of lines will envelope a caustic surface (\textit {Brennfl\"ache}),  in this case  a surface of the fourth order and class. These were precisely the types of systems of lines that led Kummer to his theory of quartic surfaces as  \textit {Brennfl\"achen} enveloped by such lines, only here the nodal line $g$ forms a double line on the quartic surface. This double line 
contains four point singularities which are
pinch points where the leaves of the surface join. They also demarcate the boundaries between real and imaginary portions of the surface as illustrated in the figure below.

\begin{figure}[ht]
        \centering 
        \includegraphics[width= 0.33 \textwidth]{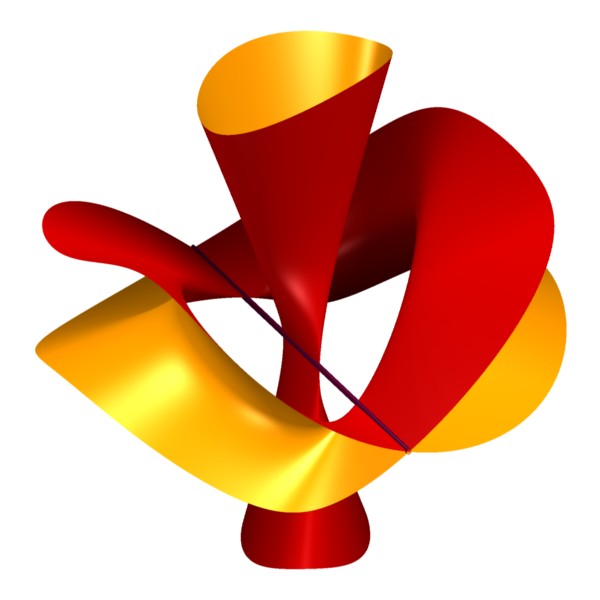}
        \caption{A Pl\"ucker complex surface showing four pinch points on its double line. Graphic courtesy of Oliver Labs.}
\label{fig:sample-image}
\end{figure}

Pl\"ucker  designed  roughly thirty-five zinc models to illustrate the many different shapes formed by these special quartics with a double line.  These models were  shown both in Germany and abroad, and in 1866 Pl\"ucker took some with him to Nottingham, where he spoke about these complex surfaces at a conference.  Arthur Cayley and other leading British mathematicians quickly took an interest in them, as they were among the more  
 exotic geometrical objects studied at this time \cite{Cayley 1871}.

Pl\"ucker also studied the local structure of  line complexes by fixing points or planes. 
 In a quadratic line complex, the lines  through a point will, in general, form a non--degenerate cone; likewise, those in a generic plane will envelope a non-degenerate conic. The exceptional points and planes are of great importance, however. If the cone through a point $P$ collapses into lines that lie in two planes, or if  the lines in a plane $\Pi$ degenerate into two point pencils, then $P$ and $\Pi$ are singular, and the locus of all such points and planes determines the {\it singularity surface} of the complex, which  in general  will be a Kummer surface. Thus, the classification of quadratic complexes turns out to be intimately tied to these special types of quartics, which were  were much studied objects during the latter third of the nineteenth century.\footnote{See, for example, the classic work on Kummer surfaces, \cite{Hudson  1905}, and the references cited therein.}
Since Pl\"ucker's {\it Neue Geometrie des Raumes}  contains no references to other literature, it would seem impossible to know whether he was aware of Kummer's concurrent papers. Klein's first publication on line geomtry, \cite{Klein1870}, however, drew all the current literature together and planted the study of Kummer surfaces firmly in the soil of quadratic line complexes.

Almost immediately after Pl\"ucker's death in 1868,  Klein became interested in the connection between the Pl\"ucker  and Kummer  surfaces, 
quartics with and without a double line $g$.\footnote{For a more detailed account of the relationship between the quartics of Pl\"ucker and Kummer, see \cite{Hudson 1905}, chapter VI.} 
Beginning with a Kummer surface, he considered three different types of Pl\"ucker quartics based on the relation of $g$ to the quadratic line complex. For $g\notin K_2$, the complex surface has eight double points and four pinch points, whereas for $g\in K_2$ there will be only four double points. Finally, if $g$ is a singular line in $K_2$, then the Pl\"ucker surface will have just two double points.
It seemed evident to Klein that one should start with a  Kummer surface, but then carry through deformations that systematically lowered the number of double points. Thus, the double line $g$ on a Pl\"ucker surface in effect absorbed eight of the 16 double points on a Kummer surface. Only  eight then remain as the others  pass over into the four pinch points on $g$.

\section*{Noether's Mapping and Tetrahedral Complexes}

During the summer of 1869, Klein got to know Max Noether, who had come to G\"ottingen to work with Clebsch. Noether had just discovered a mapping that sends the 3-parameter family of lines in a first-degree line complex to points in 3-space. Before he published this, though, Klein showed him how a similar mapping can be constructed for quadratic line complexes (\cite{Klein1921-23}, I: 89). In both cases, however, exceptional elements arise: one of the lines $\ell\in K_1$ will map to a curve rather than a point. In the case of a linear complex, this curve will be a conic, whereas for a quadratic complex a 5th-degree curve arises. Noether's mapping fascinated Klein, who wrote a lengthy paper analyzing the various types of surfaces that correspond to subfamilies of lines in a linear complex \cite{Klein1869}. Klein knew that admission to the Berlin seminar required candidates to submit a manuscript demonstrating their ability to undertake independent mathematical research, so he clearly wanted to make a good impression on Kummer and Weierstrass, the two directors of the seminar (\cite{Biermann 1988}, 279-281). 
When he met Lie soon thereafter, it must have come as a great surprise that the Norwegian, too, had come upon this type of mapping, though he interpreted this as a reciprocal relationship between two different line complexes.

In his paper, Klein noted how certain algebraic subfamilies contained  in the linear complex $K_1$ with distinguished line $\ell\in K_1$ will map into 3-space, where $C_2$ is the distinguished conic. Thus, the  lines in  $K_1$ that meet  $\ell$ will map to points that lie on $C_2$. Those lines that intersect $\ell$ in  a fixed point $P$ will map to planes tangent to the conic $C_2$, whereas the line $\ell$ itself goes over to the plane  containing $C_2$.
Klein's  idea was to consider the images of ``configurations'' in $K_1$  -- that is, 1-parameter families of lines given by an algebraic equation in line coordinates. An $n$th-order configuration will map onto a curve of the same degree in 3-space with singularities closely related to those of the configurations, which are special types of  ruled surfaces. The latter had been studied, independent of line geometry, by Arthur Cayley in \cite{Cayley 1864}. This paper  was cited by Klein along with two other more 
 recent studies  of key importance: these were  classifications of   
 4th- and 5th-degree ruled surfaces in 3-space made by Luigi Cremona in \cite{Cremona 1868} and by H. A. Schwarz in \cite{Schwarz 1867}. 
By drawing on these works, Klein aimed to  show which  ruled surfaces can arise  as configurations in a first-degree line complex. After describing his methodology, which required separate analyses of various special cases, Klein summarized his results in three tables. We will consider only the first of these, reproduced below, which deals with configurations that lie on first-degree congruences in $K_1$. 

Such congruences can be obtained simply by 
 intersecting $K_1$  with another linear complex, whereas the intersection of three linear complexes produces  the generators of a quadric surface. This 1-parameter set of lines is what Klein calls a configuration of the second degree, whereas his  investigation dealt with higher-degree  configurations,  utilizing methods introduced by  Cayley for studying skew surfaces. Congruences have two distinguished skew lines $R_1,R_2$, called directrices, with the property that all the lines of the congruence meet these two. Since this property holds  for any configuration that lies in a given congruence, this was the simplest case that Klein needed to analyze. 
Under the Noether mapping, a general linear congruence goes over to a quadric surface $S_2$, where the distinguished conic  $C_2\subset S_2$. If the distinguished line $\ell$ happens to belong  to the congruence,  then its image is the  plane $\Pi$, where  $C_2 \subset \Pi$. 
Configurations of degree $n$ will map to plane curves $C_n$, and when they belong to a congruence the curve $C_n$ 
 meets the conic $C_2$ in two fixed points, $P_1, P_2$. These two points are determined by the directrices of the congruence: $P_1=R_1\cap \Pi$ and $P_2=R_2\cap \Pi$. Klein then notes that there is no loss in generality in assuming that $\ell$ belongs to the congruence, since otherwise one can use stereographic projection to map $S_2$ onto $\Pi$.

 In Table I, transcribed below, he summarized his findings for a general case $A$ as well as for a limiting case $B$. The former corresponds to a congruence with distinct directrices, $R_1,R_2$, whereas the latter deals with a congruence for which $R_1 =R_2=R$. The upper indices for $R_1,R_2$ indicate the number of times that the curve $C_n$  passes through $P_1$ and $P_2$; hence the sum is in all cases $n$. The columns $\Sigma$ and $S$ refer to the double lines of the configuration  given by the two directrices $R_1,R_2$, counting each with its associated multiplicity. 
The columns on the left give the Cayley data associated with each of the possible cases. Thus, $p$ is the Riemannian genus of the surface, which can be obtained simply by taking plane sections, since the genus of these curves is independent of the plane chosen. The letters $n, \wp, \sigma$ stand for the degree, followed by the number of double lines, respectively, stationary lines, those where the tangent plane has 3-fold contact with the surface.

\begin{figure}[ht]
        \centering 
	   \includegraphics[width=11.5cm,height=8.5cm]{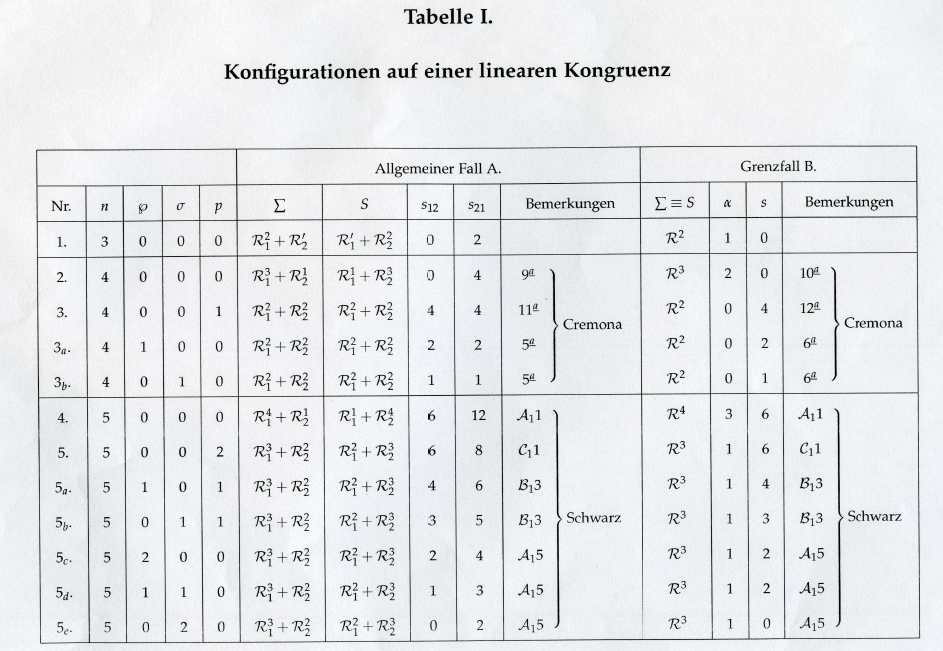}
        \caption{This table from Klein's unpublished manuscript \cite{Klein1869} summarizes his results for various configurations related to the ruled surfaces investigated by Cremona and Schwarz.}
\label{fig:sample-image}
\end{figure}

One sees at a glance from Klein's Table I that, for the general case, three types of quartic surfaces arise within Cremona's classification system, whereas four types of quintics can be obtained following Schwarz's numbering. When there is only one directrix, however, the quintic configurations are the same as before, whereas three different types of quartics emerge.

This study led to many complications, which may account for why Klein never published these results. In fact, he seems to have dropped the idea of trying to study the interior geometry of linear complexes altogether after he met Lie, who had already made considerable progress with  so-called  tetrahedral complexes. These are special quadratic line complexes $T$, which take their name from the fact that their singularity surface $S$ is the degenerate quartic given by four planes, thus a tetrahedron. More precisely, this is the point locus of $S$; as a locus of planes the locus consists of the four vertices. A generic line $\ell\in T$ 
will  intersect the faces of the singular tetrahedron in a fixed cross ratio. By varying the cross ratio, one then obtains a 1-parameter family of such complexes associated with a given tetrahedron. These complexes thus depend on 13 parameters, as opposed to 19 for a generic quadratic line complex. 

Lie typically chose this tetrahedron to consist of the three coordinate planes and the plane at infinity. Since this figure remains fixed under the 3-parameter family of transformations: $$x'=\alpha x, \qquad y'=\beta y, \qquad  z'=\gamma z,$$ he could generate the complex simply by letting this group act on a single line owing to the fact that the cross ratio will remain invariant. 
By making use of these and other mappings, Lie was able to uncover numerous  types of curves and surfaces enveloped by the lines in the complex. These geometric objects fall into distinct classes that can easily be mapped to those in another class (on this and related work of Lie and Klein, see \cite{Rowe 1989}, 230-239).

 Klein's  motivation was much the same as Lie's: both of them hoped to exploit ``transfer principles''  that would enable them to discover new  results by means of known properties of objects in 3-space. 
Lie took this approach in 
 studying the geometry of tetrahedral complexes, whereas 
 Klein  wanted to use the Noether mapping to explore the geometry of a first-degree line complex. Lie had also found this mapping, but his version of it set up a reciprocal relationship between the linear complex $K_1$ and a special quadratic complex $K_2$, namely the lines in the range space that meet a fixed conic $C_2$, the singular curve under the Noether mapping. Thus a line $\ell \in K_1$ will, in general, map to a point $P$, which together with $C_2$ forms a fixed cone in the range space. (If $P\in C_2$ then this cone collapses to the pencil of lines through $P$ in the plane of $C_2$.) Conversely, taking a line $\ell '\in K_2$, the preimages of the points $Q\in \ell '$ will be a pencil of lines in $K_1$. 

Lie's interpretation of the Noether map was only one of several innovations in line geometry that he was exploring at this time. Klein found these ideas much more promising than the quite restrictive approach he had followed, which involved comparing the preimages of special curves that corresponded to the various types of quartic and quintic ruled surfaces studied by Cremona and Schwarz. This helps to explain why no direct trace of \cite{Klein1869}, his investigation using the Noether mapping, can be found in Klein's published work. Interestingly enough, however, Lie's major paper in volume 5 of {\it Mathematische Annalen} \cite{Lie 1872} does allude to certain results on quartic surfaces that Klein found by means of the Noether mapping. This paper, a revision of Lie's dissertation, was sent to Klein in November 1871. Since the passage in question does not appear in his dissertation, this suggests the possibility that it was added by Klein, who edited several of Lie's early papers. Indeed, that was what happened, but Lie also made several other more significant changes in the text in reaction to Klein's remarks about his dissertation. Some of these matters will be taken up below, but  first let us return to the events connected with Lie's first  major breakthrough.

\section*{Lie's Line-to-Sphere Mapping}

After spending the winter semester of 1869-70 in Berlin, Lie and Klein  parted ways for roughly two months. Lie left Berlin at the end of February to spend some time with Clebsch in G\"ottingen.\footnote{On Lie's various experiences after leaving Berlin, see (\cite{Stubhaug 2002}, 137-146).} This was a few weeks before the regular semester ended, so Klein stayed on for a few more weeks. Weierstrass asked him to lecture in his seminar on Cayley's projective metric, a topic that eventually led to his famous papers on non-Euclidean geometry (\cite{Klein1871}, \cite{Klein1873}). He later joined up with Lie in Paris in late April, and soon thereafter they both got to meet Gaston Darboux, Camille Jordan, Michel Chasles, and several other leading Parisian mathematicians. 

By far the most important of these new connections was the one they made with Darboux. Despite the political tensions that followed the Franco-Prussian War and the annexation of territories in Alsace and Lorraine by the Germans, Darboux and Klein remained on friendly terms. A  letter from Klein, written on 1 November 1871, suggests why they saw one another as mathematical allies: 

\begin{quote}

Sie schrieben mir gelegentlich einmal, die neuere Geometrie erfreue sich in Frankreich noch lange nicht der gebuehrenden Anerkennung. Glauben Sie nur nicht, dass das bei uns so sehr viel anders ist. \dots 
Waehrend Steiner ein productives Genie war, dass sich gewiss nicht auf einseitige Methoden beschraenken konnte (obwohl er so that, als wenn er Analysis etc. verachtete) greifen seine Nachfolger gerade die Einseitigkeit desselben auf und verschliessen sich fuer Alles, was neben die Steiner'schen Principien von 1832 hinausgeht. – meine Arbeit ueber die Nicht-Euklidische Geometrie ist gerade aus dem Widerstreit gegen eine verkehrte Auffassung der neueren Geometrie hervorgegangen. Ich aeusserte gelegentlich, dass moeglicherweise Cayley und Lobatschefsky identisch seien; ich wurde ausgelacht. Ich hielt den Gedanken fest und schrieb ihn einigen Leuten, die mich auch eines Besseren belehrten. Das aergerte mich, und da habe ich gezeigt, dass meine Vermuthung die richtige war und die ganze Nicht-Euklidische Geschichte viel durchsichtiger macht, als dieselbe frueher war (quoted from \cite{Tobies 2016}, 110).

\end{quote}

It was during his stay in Paris that 
 Lie  gradually developed the ideas that led to his famous line-to-sphere mapping, a contact transformation with many interesting properties. Just when and how this came about has always been shrounded in mystery. Up until Klein's arrival, Lie had almost no interactions with Parisian mathematicians. Neither he nor Klein spoke French very well, but the latter was far more gregarious.
Klein 
 had a natural charm in social situations, whereas Lie was at times awkward in the company of 
 strangers. No doubt he was often lost in thought, and 
 even Klein, who spent hours with him virtually every day when they were in Paris together, experienced great difficulties when trying to follow what Lie was thinking. This much, though, seems clear. Lie spent many days contemplating the special situation that arises when the conic $C_2$ in the Noether mapping  is  displaced and becomes the ``Kugelkreis'' $C_2^*$, the imaginary circle at infinity  that lies on all spheres. In this case, Lie could interpret the images of lines in $K_1$ as spheres of radius zero, each of which forms a cone of isotropic lines determined by the fixed circle $C_2^*$. This latter setting was a familiar one to French geometers, in particular Chasles and Darboux; during the 1860s, the latter had begun to cultivate a new geometry of space based on the properties of spheres. He thus realized that what Lie had discovered was a surprising new connection between this French-style sphere geometry and Pl\"ucker's line geometry (see \cite{Darboux 1899}). 

In all likelihood, Lie was guided by the following analogy between surfaces generated by lines and those enveloped by spheres. In the context of line geometry, the intersection of three first degree line complexes will be a one-parameter family of lines that forms a system of generators for a quadric surface, for example a hyperboloid of one sheet in the real case. The analogous objects in sphere geometry are the special quartic surfaces known as Dupin cyclides, and indeed Lie's mapping will transform each of these two types of surfaces to the other.  To see this most easily, one can take three skew lines in space that  map to three spheres. Then the set of  lines that meet these mutually skew lines forms a system of generators of a quadric surface, and since Lie's mapping is a contact transformation these lines go over to a one-parameter family of spheres tangent to three fixed spheres. The second system of generators will have the same property, and the analogy is then clear, since a Dupin cyclide is enveloped by two families of spheres. 
 Moreover, this mapping has the property that the asymptotic curves of the first surface go over to the curvature lines of the second. In this case, the generators themselves are the asymptotic curves, and these then correspond to the circles of tangency of the Dupin surface (see \cite{Lie and Scheffers 1896}, pp. 470-475).

Klein and Lie soon learned that the French geometers had come up with a far-reaching extension of these older findings that led to a theory of surfaces based on spheres.
In 1864 Darboux and Theodore Moutard began work on generalized cyclides, which they studied in the context of inversive geometry. Klein and Lie learned about this new French theory when they met Darboux just before the Franco-Prussian War broke out. Darboux later developed the theory of generalized cyclides by introducing pentaspherical coordinates.\footnote{Darboux had already worked out many of these ideas when Klein and Lie met him, but it took him another two years to develop the whole theory in detail and publish it in \cite{Darboux 1873}. For further background on his early career and mathematical research, see \cite{Croizat 2016}.}
 These cyclides are special quartic surfaces with the property that they meet the plane at infinity in a double curve, namely the imaginary circle $C_2^*$ that lies on all spheres. Darboux also found that their lines of curvature are algebraic curves of degree eight. This finding set up one of Lie's earliest discoveries, communicated to Darboux at that time. 

This breakthrough came from Lie's line-to-sphere mapping, when he considered the caustic surface enveloped by lines in a congruence of the second order and class. These were, of course,  Kummer surfaces; indeed, this was the original context in which the Berlin mathematician had studied them. Lie realized that his mapping established a natural reciprocity between geometrical objects in two spaces $E_1$ and $E_2$ in which the two line complexes, $K_1$ and $K_2$, resided. Thus, a surface $F_2\subset E_2$, viewed as a point locus, will pull back to a 2-parameter family of lines, thus a congruence in $E_1$, and vice-versa. The caustic surface $F_1\subset E_1$ corresponding to a generalized cyclide $F_2$ turns out to be a  Kummer surface. Since 
Darboux had shown how to determine the line of curvature on cyclides, Lie  immediately deduced that the asymptotic curves on a Kummer surface could be found and that these are algebraic of degree sixteen. 

Lie and Klein discussed this breakthrough in detail, as Klein gradually came to understand Lie's line-to-sphere mapping. He
 quickly realized that Lie's claim regarding the asymptotic curves of Kummer surfaces was  correct; in fact,  
he had already come across these same curves  in his own work. Klein had found them while studying quadratic line complexes that share the same Kummer surface as their surface of singularity, a paper he completed in June 1869 \cite{Klein1870}. At that time, however,  he had not realized  that these were the asymptotic curves on the surface. 
Now Klein was able to confirm not only that Lie's general conclusion was correct but that he could deduce other properties of these curves by exploiting his theory of one-parameter families of quadratic line complexes with a fixed singularity surface. Klein studied these, in analogy with confocal systems of surfaces, by means of 
 elliptic coordinates, which had been introduced by Jacobi in connection with surface theory. By varying the parameter associated with the quadratic complexes, Klein was able to  sweep out all the asymptotic curves that lie on a fixed Kummer surface.

\begin{figure}[h]
        \centering 
        \includegraphics[width=9.5cm,height=6.5cm]{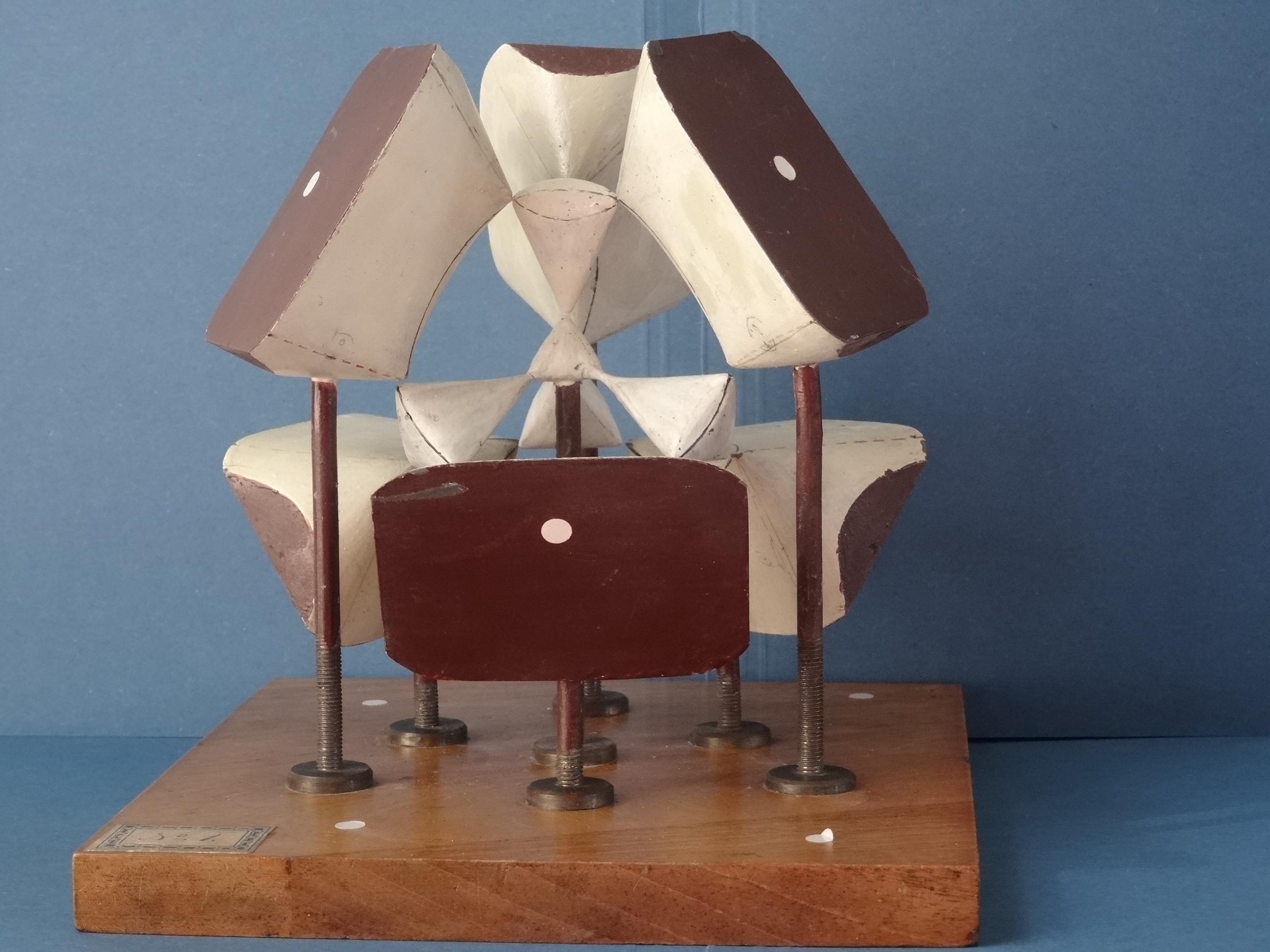}
        \caption{Klein's model of a Kummer surface designed in 1871. Courtesy of the Collection of Mathematical Models,  G\"ottingen University.}
\label{fig:maya}
\end{figure}

Klein and Lie only got to discuss these matters for about two weeks, however, because in mid-July 
France declared war on Prussia.  Klein hastily fled from Paris, but from back home in D\"usseldorf  he reported to Lie that he was able to trace the paths of these asymptotic curves and  to describe their singularities. He managed to visualize these with the help of  a physical model 
of a Kummer surface made by his friend Albert Wenker. By this time, Klein realized that what he had told Lie earlier in Paris about the singularities of these asymptotic curves was, in fact, incorrect. After giving the necessary corrections in a letter from 29 July, he added:

\begin{quote}

I came across these things by means of Wenker's model, on which I
wanted to sketch asymptotic curves. To give you a sort of intuitive idea
how such curves look, I enclose a sketch. The Kummer surface contains
hyperboloid parts, like those sketched; these are bounded by two of the
six conics ($K_1$ and $K_2$) and extend from one double point ($d_1$) to another
($d_2$). Two of the curves are drawn more boldly; these are the two that not
only belong to linear complexes but also are curves with four-point contact.
They pass through $d_1$  and $d_2$  readily, whereas the remaining curves have
cusps there. This is also evident from the model. At the same time, one sees
how $K_1$  and $K_2$  are true enveloping curves.\footnote{Klein to Lie, 29 July 1870, {\it Letters from Felix Klein to Sophus Lie, 1870-1877}, Heidelberg:  Springer-Verlag, scheduled to appear in 2018.}
\end{quote}

\begin{figure}[h]
        \centering 
        \includegraphics[width=6.5cm,height=9.5cm]{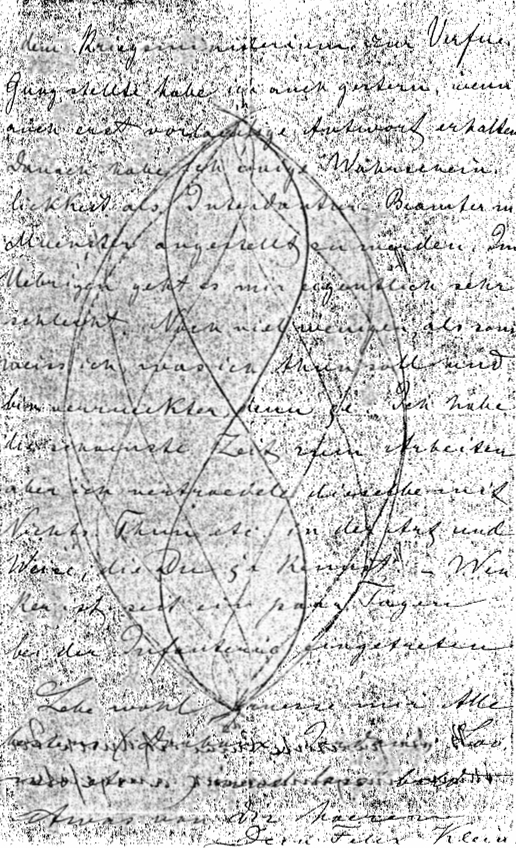}
        \caption{Klein's sketch of the asymptotic curves between two double points on a Kummer surface from his letter to Sophus Lie.}
\label{fig:maya}
\end{figure}

By the ``hyperboloid parts'' on a Kummer surface, Klein meant those places where the curvature was negative. Only in these regions were the asymptotic curves real and hence visible. He later reproduced the  same figure in the note that he and Lie sent to Kummer for publication in the {\it Monatsberichte} of the Prussian Academy \cite{Klein and Lie 1870}. 
As briefly noted above, Klein later designed four models in 1871, one representing a general Kummer surface with 16 real nodes, and three others for Pl\"ucker surfaces with eight, four, and two singular points, respectively. His model with 16 nodes is shown in figure 6.

\section*{Lie's Applications of the Line-to-Sphere Mapping}

Already in early July 1870,  thus before Klein left Paris, Lie communicated some of his findings to the Norwegian Scientific Society in Christiania, though this note was only published by Ludwig Sylow in 1899, the year of Lie's death (\cite {Lie 1934},  86-87). In it Lie simply announced a series of results, one of which relates to ruled surfaces, i.e. configurations, that lie in a linear complex, the topic of Klein's unpublished manuscript \cite{Klein1869}, discussed above. The asymptotic curves of such surfaces, Lie wrote, can all be determined by quadrature. A short time later, Lie wrote a somewhat more extensive note that Chasles submitted to the Paris Academy for publication in the {\it Comptes Rendus} on 31 October 1870 (\cite {Lie 1934}, 88-92). Here he briefly explains how he found these asymptotic curves  by using his line-to-sphere mapping; he also noted its close relationship with the Noether map, which Klein had studied just before he and Lie met in Berlin. 

Lie's result exploited another finding due to Darboux, which pertains to curvature lines on an arbitrary (smooth) surface $S$. Darboux showed that the developable surfaces $D(S)$ obtaining by taking the lines tangent to $S$ that meet $C_2^*$, the circle at infinity, will touch $S$ along an (imaginary) line of curvature. Thus, by pulling back this surface $S$, Lie observed that the caustic surface determined by the resulting congruence has one known asymptotic curve, whose tangents belong to a linear complex. Moreover, this will be true for any configuration that lies in a linear complex. Finally, Lie noted that Bonnet and Clebsch had shown that from one asymptotic curve on a ruled surface one could derive all the others. So by exploiting his mapping and these recent results obtained by other geometers, Lie had found a remarkable, though rather theoretical result. Still, all along he was thinking a good deal about how to actually integrate the curves and surfaces associated with a line complex.

Around this same time, in late October 1870, Lie submitted yet another short note to the Norwegian Academy (\cite {Lie 1934}, 93-96). His central theme here concerned the connections between line complexes and first-order partial differential equations. These ideas were of central importance for Lie (see \cite{Hawkins 2000}, 20-42), 
 but only toward the end of his life did he elaborate on them with the help of Georg Scheffers \cite{Lie and Scheffers 1896}. Nevertheless, this note already contains a synopsis of some of the major results that would occupy Lie's attention throughout the next year. First, he noted that  partial differential equations that define a line complex are those for which the characteristic curves are asymptotic curves of the integral surfaces. Second, he discussed the differential equations that determine lines of curvature and geodesic curves, pointing out that solving any one of these three types of equations leads to solutions of the others by means of appropriate transformations. Finally, he alluded to Jacobi's determination of the geodesic curves on quadric surfaces as a concrete example of the theory he had in mind. 

Lie elaborated on this final idea in the second half of the text he planned to submit as his doctoral dissertation the following year. His doctoral committee, however, insisted on following the strict rules, which stipulated that candidates for a doctorate submit their work either in Latin or Norwegian. So Lie translated the first half of this text  into his mother tongue, and defended his thesis in June 1871. According to Elling Holst, the two mathematicians who were tasked to report on Lie's work could barely understand a word of it (\cite {Stubhaug 2002}, 156-157). Klein and Clebsch were, of course, eager to publish a revised version in {\it Mathematische Annalen}. When  Klein received the second (untranslated) part of Lie's German text, he noted that the ruled quartics Lie had found as special singularity surfaces  could be identified using Cremona's classification. This part of Lie's paper, in fact, was originally very sketchy (see \cite {Lie 1934}, pp. 194-199), but Klein induced Lie to give a far more elaborate account of the various types of quadratic line complexes that fell under his new theory. These had ruled quartics with two double lines as their 
 singularity surfaces, and they were distinguished by the property that the asymptotic curves of these quartics could be transformed to the geodesics on a quadric surface. The latter, following Jacobi, could be found by means of hyperelliptic integrals, the reduction Lie had announced in his note for the Norwegian Academy. 

Lie's revised dissertation \cite{Lie 1872} was arguably the most brilliant work he ever wrote. Unfortunately,  very few who tried to read it  got very far. Lie was the sort of mathematician who spent more time finding new results than explaining what he had already found to others, and of course his 
 ideas  were far ahead of their time. A glance at the  fourth and final part of \cite{Lie 1872} shows that Klein's ideas and input played an important part in the revision Lie undertook. In fact,  Lie was very generous in praising him from the very beginning, thanking Klein for numerous ideas, more than he could possibly cite in the text (\cite {Lie 1934}, p. 3). 

In \S 23 Lie worked through 26 cases of special quadratic line complexes, noting that this led to 24 independent types in the end. He also identified the various types of ruled surfaces that were associated with them. This elaboration was carried out in response to a letter from Klein, written on 29 July 1871, in which the latter wrote about two questions connected with Lie's theory: 
``a) To look at what kind of detailed results the general reasoning about the
Kummer surface yields for these [special singularity] surfaces.
b) To give all such surfaces (and complexes of the second degree).''\footnote{Klein to Lie, 29 July 1871, {\it Letters from Felix Klein to Sophus Lie, 1870-1877}, Heidelberg:  Springer-Verlag, scheduled to appear in 2018.}

The second question was one Klein had touched upon in his dissertation, written in Bonn  in 1868  when he was still only nineteen years old (\cite{Rowe 1989}, 217-220). His mentor then was  Rudolf Lipschitz, owing to the fact that Pl\"ucker had already died, and it was he who prompted Klein to use Weierstrass's new theory of elementary divisors in his work. Klein clearly recognized that this was an important tool for 
 the classification of quadratic lines complexes, a problem that  remained very much on his mind. Beyond the purely algebraic aspect, however, he was even more interested in the geometric problem that had originally inspired Pl\"ucker. So the question he posed to Lie was related to the issue of 
 visualizing the various types of quadratic line complexes within a given classification scheme. 

Klein gave this problem to his student Adolf Weiler, whose dissertation from 1873 was published in \cite{Weiler 1874}. Weiler found 48 distinct types of degenerate quadratic line complexes based on the geometry of their singularity surfaces. Only in the case of distinct eigenvalues will this quartic be a Kummer surface, in which case the complex depends on 19 parameters. The simplest degeneration leads to a Pl\"ucker  surface, which belongs to a complex with 18 parameters. Among Cremona's ruled quartics,
five of the twelve types can arise as  singularity surfaces, but not all of these are Lie quartics. Lie's  most general case belongs to a complex with 16 parameters, whereas Cremona's general type XI is the singularity surface for a complex that depends on 17 parameters. Another special quartic listed by Weiler arises in connection with a complex with 16 parameters, namely the famous Roman surface of Steiner that Kummer found independently, as described above. 

In the table summarizing these results, Weiler lists the cases enumerated by Lie in \S 23 of \cite{Lie 1872} alongside the corresponding numbers in Cremona's general classification of ruled quartics. One can easily see how this  reflects not only Klein's desire to call attention to Lie's results but also his larger interest in this topic, going back even to his work on the Noether mapping, which brought a different set of Cremona quartics to light. 

All of these tightly bound threads  connected with  line geometry took place over a very short time, but  by the mid-1870s both Klein and Lie had moved on  to other areas of research. Still, they never lost sight of the work they had done together in their youth. In fact, their divergent views about their respective roles during that time became a major factor underlying  the dramatic break in their friendship that occurred in 1893. Afterward, Lie struck up a 
 collaboration with Scheffers that  represents his first sustained attempt to look backward in order to elaborate on the germinal ideas that eventually led him to his work on transformation groups and differential equations. Klein, on the other hand, had made various attempts to showcase his early work with Lie already some years before. His motivation was sparked at least in part by the work of young Italian geometers, in particular Corrado Segre.
It was hardly an accident that Segre played a major role in bringing the text of Klein's largely forgotten ``Erlanger Programm'' back to light, but that is part  of another story (see \cite{Hawkins 1984}, \cite{Luciano and Roero 2012},  and  \cite{Rowe 2017a}).

\section*{Appendix}

Lie Papers, Oslo,  Brevs. Nr. 289

\bigskip

Clebsch an Lie, G\"ottingen, den 30. Oktober 1869
\bigskip
\noindent

Geehrter Herr Lie!
\bigskip

Die beiden S\"atzen, welche Sie so freundlich waren, mir mitzutheilen, sind mir nicht bekannt, insbesondere der zweite nicht; f\"ur diesen konnte ich mir einen kurzen und eleganten analytischen Beweis machen, welcher es gestattet von der Darstellung der Steinerschen Fl\"ache durch rationaler Funktionen zweier Parameter zu der analogen Darstellung der zweiten Steinerschen Fl\"ache \"uberzugehen. 

Was den ersten Satz betrifft, so hat sich Dr. Klein (welcher das Pl\"uckersche Werk fortgesetzt hat) voriges Sommer hier mit solchen Sachen besch\"aftigt, und ist m\"oglicherweise zu gleichen Resultaten gekommen. Ich bezweifele nicht, dass die Fl\"ache allgemein sein wird; vielleicht theilen Sie mir gelegentlich n\"aheres mit. Doch w\"urde ich Ihnen vor allem empfehlen Dr. Klein, der jetzt in Berlin ist (Carlstrasse 11) aufzusuchen und bestens von mir zu gr\"ussen; Sie werden in ihm einen strebsamen und liebensw\"urdigen Mann finden. Auch kennt er Sie schon durch Ihre Abhandlung, f\"ur diese freundliche \"Ubersendung ich Ihnen bei dieser Gelegenheit meinen herzlichen Dank sage.

Dr. Klein wird Ihnen auch am besten sagen k\"onnen, wie es mit der Pl\"uckerschen Stelle p. 222 ist. Ich w\"usste nicht, was Ihrer Schlussweise entgegenzusetzen w\"are, aber ich bin augenblicklich in diesen Sachen nicht so zu Hause, und kann mich im Augenblick auch nicht so hineinarbeiten, um sicher entscheiden zu k\"onnen. Mit Dr. Klein w\"urden Sie sich gleich 
verst\"andigen. 
\bigskip

In der Hoffnung weiter von Ihnen zu h\"oren.
\bigskip

Ihr ergebener
\bigskip

A. Clebsch

\end{document}